\documentclass[12pt]{article}

\usepackage{amssymb}
\usepackage{amsmath}
\usepackage{amsthm}

\begin{document}

\title{Irreducible Characters for Algebraic Groups in
Characteristic Three}

\date{}
\author{Zhongguo Zhou\quad  Xiangqin Meng}
\date{\small\it College of Science, Hohai University\\
 Nanjing,  210098, P.R.China.\\ {
{\rm e-mail: zhgzhou@hhu.edu.cn, mxqzhg@hhu.edu.cn}}} \maketitle

\begin{abstract}
In this note, we determine the irreducible
characters for the simple algebraic groups of type $A_5$  over
an algebraically closed field $K$ of characteristic $3$, by using a
theorem of Xi Nanhua \cite{X1} and the Matlab software. In order to
obtain higher speed than in \cite{YZ2001,YZ2002} we modify
 the algorithm to compute the irreducible characters.
\end{abstract}





The determination of all irreducible characters is a big theme in
the modular representations of algebraic groups and related finite
groups of Lie type. But so far only a little is known concerning
it in the case when the characteristic of the base field is less
than  the Coxeter number.

Gilkey-Seitz gave an algorithm to compute part of characters of
$L(\lambda)$'s with $\lambda \in X_1(T)$ for $G$ being of type $G_2$,
$F_4$, $E_6$, $E_7$ and $E_8$ in characteristic $2$ and even in
larger primes in \cite{gs}. Dowd and Sin gave all characters of
$L(\lambda)$'s with $\lambda \in X_1(T)$ for all groups of rank less than
or equal to $4$ in characteristic $2$ in \cite{DS}. They got their
results by using the \lq\lq standard\rq\rq ~Gilkey-Seitz algorithm
and computer.  L. Scott et al.
computes the characters for $A_4$ when $p=5,p=7 $
by computing  the maximal submodule in a baby Verma module \cite{sc}.
Anders Buch and Niels Lauritzen
also obtain this result for $A_4$ when $p=5 $ with Jantzen's sum formula \cite{bs}.

An element ${\mathfrak  x}_{(p^n-1)\rho-\lambda} \in {\mathfrak u}^-_n$ for
each irreducible module $L(\lambda)$ with $\lambda \in X_n(T)$ was defined
in [5, \S 39.1, p. 304] and  [7, p. 239]. This element could be
used in constructing a certain basis for $L(\lambda)$, computing $\dim
L(\lambda)$, and determining ${\text{\rm ch}}(L(\lambda))$. In this way, Xu and Ye, Ye and Zhou determined
all irreducible characters for the special linear groups
$SL(5,K)$, $SL(6,K)$ and $SL(7, K)$, the special orthogonal group
$SO(7, K)$ and the symplectic group $Sp(6, K)$ over an
algebraically closed field $K$ of characteristic $2$ in \cite{XY1997,YZ2000} and for the special orthogonal group $SO(7, K)$ and the
symplectic group $Sp(6, K)$ over an algebraically closed field $K$
of characteristic $3$ in \cite{YZ2001,YZ2002}. However, it need so much time to compute the irreducible
characters when the characteristic is bigger than 3 that this become to be impossible mission.
So we must find out faster algorithm to finish this work.  In the present note, we  shall work
out all irreducible characters for the simple algebraic groups of
type $A_5$  over an algebraically closed field $K$ of
characteristic $3$ and introduce  how to obtain faster speed with modified algorithm.
Of course we will explain why
our results  should be right. We shall freely use the notations in [9]
without further comments.

\section{PRELIMINARIES}
\label{}

Let $G$ be the simple algebraic group of type $A_5$  over an algebraically
closed field $K$ of characteristic $3$. Take a Borel subgroup $B$ and a maximal torus
$T$ of $G$ with $T\subset B$. Let $X(T)$ be the character group of $T$, which is also
called the weight lattice of $G$ with respect to $T$. Let $R\subset X(T)$ be the root
system associated to $(G,T)$, and choose a positive root system $R_+$ in such a way that
$-R_+$ corresponds to $B$. Let
$$
S=\{\alpha_1, \alpha_2, \alpha_3, \alpha_4,\alpha_5 \}
$$
be the set of simple roots of $G$ such that
$$
R_+=\{\alpha_1, \alpha_2, \alpha_3, \alpha_4,\alpha_5, \alpha_{ij}=\alpha_i +\cdots+\alpha_j,\, 1\leq i<j\leq 5\}.
$$
Let $\omega_i (1 \leq i \leq 5)$ be the fundamental weights of $G$ such
that $\langle \omega_i, \alpha^{\vee}_j \rangle=\delta_{i j}$, the Kronecker delta,
and denote by $\lambda=(\lambda_1, \lambda_2, \lambda_3, \lambda_4, \lambda_5)$ the weight $\lambda=
\lambda_1\omega_1 + \lambda_2\omega_2 + \lambda_3\omega_3 + \lambda_4\omega_4+ \lambda_5\omega_5$ with $\lambda_1, \lambda_2, \lambda_3,
\lambda_4 , \lambda_5\in \Bbb Z$, the integer ring. Then the dominant weight set is as
follows:
$$
X(T)_+=\{ (\lambda_1, \lambda_2, \lambda_3, \lambda_4, \lambda_5) \in X(T) \mid \lambda_1, \lambda_2,
\lambda_3, \lambda_4, \lambda_5 \geq 0 \}.
$$
\par
Let $W={N_G(T)}/T$ be the Weyl group and let $W_3$ be the affine Weyl
group of $G$. It is well-known that for $\lambda \in X(T)_+$, $H^0(\lambda)$
is the induced $G$-module from the 1-dimensional $B$-module $ K_{\lambda}$
which contains a unique irreducible $G$-submodule $L(\lambda)$ of the highest
weight $\lambda $. In this way, $X(T)_+$ parameterizes the finite-dimensional
irreducible $G$-modules. We set ${\text{\rm ch}} (\lambda)={\text{\rm ch}}(H^0(\lambda))$ and
${\text{\rm ch}_3}(\lambda)={\text{\rm ch}}(L(\lambda))$ for all $\lambda \in X(T)_+$.
Moreover, ${\text{\rm ch}} (\lambda)$ is given by the Weyl
character formula, and for $\lambda \in X(T)_+$, we have
$$
{\text{\rm ch}}(\lambda)={{\sum_{w\in W} det(w)e(w(\lambda +\rho))} \over {\sum_{w\in
W}det(w)e(w \rho )}}.
$$
For $\lambda=(a, b, c, d,e) \in X_1(T)$,  we have
$$\begin{array}{lllllll}
\dim H^0(a, b, c, d, e)=
&\frac{1}{ {2^83^35}}(a+1)(b+1)(c+1)(d+1)(e+1))(a+b+2)\\
                  &(b+c+2)(c+d+2)(d+e+2)(a+b+c+3)\\
                  &(b+c+d+3)(c+d+e+3)(a+b+c+d+4)\\
                  &(b+c+d+e+4)(a+b+c+d+e+5),
\end{array}$$
Let $F^n$ be the $ n$-th Frobenius morphism of $G$ with $G_n \subset G$ the
scheme-theoretic kernel of $F^n$. Let $V^{[n]}$ be the Frobenius twist for
any $G$-module $V$. It is well-known that $ V^{[n]}$ is trivial regarding as
a $G_n$-module. Moreover, any $G$-module $M$ has such a form if the action
of $G_n$ on $M$ is trivial. Let
$$
X_n(T)=\{ (\lambda_1, \lambda_2, \lambda_3, \lambda_4, \lambda_5) \in X(T)_+ \mid \lambda_1, \lambda_2,
\lambda_3, \lambda_4 , \lambda_5< 3^n \}.
$$
Then the irreducible $G$-modules $L(\lambda)$'s with $\lambda \in X_n(T)$ remain
irreducible regarded as the $G_n$-modules. On the other hand, any irreducible
$G_n$-module is isomorphic to exactly one of them.
\par
For $\lambda \in X(T)_+$, we have the unique decomposition
$$
\lambda=\lambda^0 + 3^n\lambda^1 \qquad {\text {\rm with}}\quad \lambda^0 \in X_n(T),\,
\lambda^1 \in X(T)_+.
$$
Then the Steinberg tensor product theorem tells us that
$$
L(\lambda) \cong L(\lambda^0) \otimes L(\lambda^1)^{[n]}.
$$
Therefore we can determine all the characters ${\text{\rm ch}_3\,}(\lambda)$ with $\lambda \in
X(T)_+$ by using the Steinberg tensor product theorem, provided that all
the characters ${\text{\rm ch}_3\,}(\lambda)$ with $\lambda \in X_1(T)$ are known.
\par
Recall the strong linkage principle in \cite{h}. We define a strong linkage
relation $\mu \uparrow \lambda $ in $X_+(T)$ if $L(\mu)$ occurs as a composition
factor in $ H^0(\lambda)$. Then $ H^0(\lambda)$  is irreducible when $\lambda $ is a
minimal weight in $X(T)_+$ with respect to the partial ordering  determined by
the strong linkage relations.
\par
Let $\mathfrak g $ be the simple Lie algebra over $\mathbb C$ which has the same
type as $G$, and $\mathfrak U$ the universal enveloping algebra of $\mathfrak g $. Let
$e_\alpha, f_\alpha, h_i (\alpha \in R_+, i=1, 2, 3,4,5)$ be a Chevalley basis of $\mathfrak g$.
We also denote $e_{\alpha_I}, f_{\alpha_I}$ by $e_I, f_I$, respectively, where
$I \in \mathcal A=\{1, 2, 3, 4,5, 12, 23, 34,45, 13, 24,35,14,25,15 \}$ The Kostant $\mathbb Z$-form $\mathfrak U_{\mathbb Z}$
of $\mathfrak U$ is the $\mathbb Z$-subalgebra of $\mathfrak U$ generated by the elements $e_{\alpha}^{(k)}:
=e_{\alpha}^k/k!, f_{\alpha}^{(k)}:=f_{\alpha}^k/k!$ for $\alpha \in R_+$ and $k \in
\mathbb Z_+$. Set
$$
{{h_i+c}\choose k}:=\frac{(h_i+c)(h_i+c-1)\cdots (h_i+c-k+1)}{k!}.
$$
Then ${{h_i+c}\choose k}\in \mathfrak U_{\Bbb Z},$ for $ i=1, 2, 3, 4,5$, $c\in\Bbb
Z$, $k\in\Bbb Z_+$. Define $\mathfrak U_k:=\mathfrak U_{\Bbb Z}\otimes K$ and call $\mathfrak U_k$ the
hyperalgebra over $K$ associated to $\mathfrak g$. Let $\mathfrak U_k^+, \mathfrak U_k^-, \mathfrak U_k^0$ be the
positive part, negative part, zero part of $\mathfrak U_k$, respectively. They are
generated by $e_{\alpha}^{(k)}$, $f_{\alpha}^{(k)}$ and ${{h_i} \choose k}$,
respectively. By abuse of notations, the images in $\mathfrak U_k$ of $e_{\alpha}^{(k)},
f_{\alpha}^{(k)}, {{h_i+c}\choose k}$, etc. will be denoted by the same notations,
respectively. The algebra $\mathfrak U_k$ is a Hopf algebra, and $\mathfrak U_k$ has a triangular
decomposition $\mathfrak U_k=\mathfrak U_k^-\mathfrak U_k^0\mathfrak U_k^+$. Given a positive integer $n$, let $\mathfrak U_n$
be the subalgebra of $\mathfrak U_k $ generated by the elements $ e_{\alpha}^{(k)}, f_{\alpha}^
{(k)}, {{h_i}\choose k}$ for $\alpha \in R_+, i=1, 2, 3, 4,5$ and $ 0 \leq k <
3^n$. In particular, $\mathfrak U=\mathfrak U_1$ is precisely the restricted enveloping algebra
of $\mathfrak g$.
Denote by $\mathfrak U_n^+, \mathfrak U_n^-, \mathfrak U_n^0 $ the positive part, negative part, zero
part of $\mathfrak U_n$, respectively. Then we have also a triangular decomposition
 $\mathfrak U_n=\mathfrak U_n^-\mathfrak U_n^0\mathfrak U_n^+$. Given an ordering in $ R_+$, it is known that the
PBW-type bases for $\mathfrak U_k $ resp. for $\mathfrak U_n$ have the form of
$$
\prod_{\alpha\in R_+}f_{\alpha}^{(a_\alpha)}\prod_{i=1}^5{{h_i}\choose b_i}
\prod_{\alpha\in R_+}e_{\alpha}^{(c_\alpha)}
$$
with $a_\alpha, b_i, c_\alpha \in \Bbb Z_+$ resp. with $0 \leq a_\alpha, b_i, c_\alpha <
3^n$.
\par
Let $\lambda=(\lambda_1, \lambda_2, \lambda_3, \lambda_4,\lambda_5) \in X_n(T)$. We set $\lambda_I=\sum_{i \in I}\lambda_i$
for $I \in \mathcal A,$ here each element $I$ is also viewed as a certain set of
simple roots. Following \cite{L2} and \cite{X1}, we define an elements ${\frak x}_{\lambda}$ in $\mathfrak U_n^- $ by
$$\begin{array}{lll}
{\frak x}_{\lambda}=&f_1^{(\lambda_{1})}f_2^{(\lambda_{1 2})}f_3^{(\lambda_{13})}f_4^{(\lambda_{14})}f_5^{(\lambda_{15})}
f_1^{(\lambda_{2})}f_2^{(\lambda_{2 3})}f_3^{(\lambda_{24})}f_4^{(\lambda_{25})}\\
& f_1^{(\lambda_{3})}f_2^{(\lambda_{34})}f_3^{(\lambda_{35})}
f_1^{(\lambda_4)}f_2^{(\lambda_{45})}
f_1^{(\lambda_{5})}.
\end{array}$$
As a special case of [7, Theorems 6.5 and 6.7], we have
 {Theorem 1}
Assume that $\mathfrak g$ is a simple Lie algebra of the simple algebraic group
of type  $A_5$  over an algebraically closed field $K$
of characteristic $3$. Let $\lambda=(\lambda_1, \lambda_2,  \lambda_3,  \lambda_4,\lambda_5) \in X_n(T)$.
\par
\noindent {\text {\rm (i)}} The element ${\frak x}_{\lambda}$ lies in $\mathfrak U_n^- $.
\par
\noindent {\text {\rm (ii)}} Let $\frak J_{\lambda}$ be the left ideal of $\mathfrak U_k$
generated by the elements $e_i^{(k)}, {{h_i}\choose k}-{{\langle \lambda,
\alpha_i^\vee \rangle}\choose k}, f_i^{(k_i)}$ $(i=1, 2, 3, 4,5$, $ k\geq 1$,
$ k_i\geq 3^{n})$ and the elements $f\in\mathfrak U_n^-$ with $f{\frak x}_{(3^n-1)\rho-\lambda}=0$.
Then $\mathfrak U_k/\frak J_{\lambda} \cong L(\lambda)$ (Note that $L(\lambda)$ has a $\mathfrak U_k$-module
structure, which is irreducible).
\par
\noindent {\text {\rm (iii)}} As a $\mathfrak U_n^-$-module, $L(\lambda)$ is isomorphic to
$\mathfrak U_n^-{\frak x}_{(3^n-1)\rho-\lambda}$.

By abuse of notations, the images in $\mathfrak U_k/\frak J_{\lambda} \cong L(\lambda)$ of
$f_i^{(k_i)}$ and $f_I^{(k_I)}$ will be denoted by the same notations. We
shall use this theorem to computer the multiplicities of the weight spaces
for all the dominant weight of $L(\lambda)$, to compute $\dim L(\lambda)$, and to
determine ${\text{\rm ch}}(L(\lambda))={\text{\rm ch}}_3(\lambda)$ ($\lambda \in X_1(T)$) in this note, when $G$ is the simple
algebraic group of type  $A_5$.

\section{ CHARACTERS OF THE IRREDUCIBLE MODULES OF $G$ }

From now on we shall assume that $n=1$. Denote by $V^*$ the dual
module of $V$, then we have by the duality that ${\text{\rm ch}}
H^0(\lambda)^*={\text{\rm ch}}(-w_0\lambda)$, and ${\text{\rm ch}}
L(\lambda)^*={\text{\rm ch}}_3(-w_0\lambda)$.
Furthermore, the elements $f_I$ ($I \in {\mathcal A}$ ) satisfy the
following commutator relations:
$$\begin{array}{cccccccc}
f_1f_2&=f_2f_1+f_{12},\quad \qquad \;\;\;
f_2f_3=f_3f_2+f_{23},\\
f_3f_4&=f_4f_3+f_{34},\qquad \quad \;\;
f_{12}f_3=f_3f_{12}+f_{123},\\
f_{23}f_4&=f_4f_{23}+f_{234},\qquad \quad
f_1f_{23}=f_{23}f_1+f_{123},\\
f_2f_{34}&=f_{34}f_2+f_{234},\qquad \quad \!
f_1f_{234}=f_{234}f_1+f_{1234},\\
f_{12}f_{34}&=f_{34}f_{12}+f_{1234},\qquad
f_{123}f_4=f_4f_{123}+f_{1234},\\
f_If_{I'}&=f_{I'}f_I \qquad \text {for ~~all ~~the ~~other ~~
$I, I' \in \mathcal A$}.
\end{array}$$
\par
Now we can obtain our main theorems. Let $e(\nu) = \sum_{w \in W_{\nu}}
w(\nu)$ be the sum of weights of the W-orbit of $\nu$ for all $\nu \in
X(T)_+$. It is well-known that $\{ {\text{\rm ch}}(\nu) \vert \nu \in X(T)_+ \}$,
$\{ {\text{\rm ch}}_3(\nu) \vert \nu \in X(T)_+ \}$ and $\{ e(\nu) \vert \nu \in
X(T)_+ \}$ form bases of $\mathbb Z[X(T)]^W$, the W-invariant subring of
$\mathbb Z[X(T)]$, respectively. According to the Weyl character formula and
the Freudenthal multiplicity formula, we get a change of basis matrix
$A=(a_{\lambda\nu})_{\lambda,\nu \in X(T)_+}$ from $\{ e(\nu) \vert \nu \in X(T)_+
\}$ to $\{ {\text{\rm ch}}(\nu)\vert \nu \in X(T)_+ \}$, which is a triangular matrix
with $1$ on its diagonal, i.e.
$$
{\text{\rm ch}}(\lambda)=\sum_{\nu \prec \lambda,~ \nu \in X(T)_+}a_{\lambda\nu}e(\nu)
$$
with $a_{\lambda\lambda}=1$(cf. [12]). Based on our computation, we get
another change of basis matrix $B=(b_{\lambda\nu})_{\lambda,\nu \in X(T)_+}$
\!from $\{ e(\nu) \vert \nu \!\in\! X(T)_+ \}$ to $\{ {\text{\rm ch}}_3(\nu) \vert$ \linebreak$ \nu \in
X(T)_+ \}$, which is also a triangular matrix with $1$ on its diagonal.
\par
Let us mention our computation of $B$ more detailed. First of all,
we compute ${\frak x}_{2\rho -\lambda}$ for any $\lambda \in
X_1(T)$. It is well known that for each dominant weight $\nu$ of
$H^0(\lambda)$, $\beta=\lambda-\nu$ can be expressed in terms of sum
of positive roots, and there exist many ways to do so. Each way
corresponds to an element $f_{\beta}{\frak x}_{2\rho - \lambda}$ in
$\mathfrak U_n$. Then we compute various $f_{\beta}{\frak x}_{2\rho -
\lambda}$. Note that each $f_{\beta}{\frak x}_{2\rho - \lambda}$ can
be written as a linear combination of the basis elements of $\mathfrak U_n$
with non-negative integer coefficients, and the typical images of
all non-zero $f_{\beta}{\frak x}_{2\rho - \lambda}$'s generate the
weight space $L(\lambda)_{\nu}$ of the irreducible submodule
$L(\lambda)$ of $H^0(\lambda)$. Therefore, we can easily determine
the dimension of $L(\lambda)_{\nu}$, provided that we compute the
rank of the set of all these non-zero $f_{\beta}{\frak x}_{2\rho -
\lambda}$'s. It can be reduced to compute the rank of a
corresponding matrix. Finally, we obtain the formal character of
$L(\lambda)$, which can be written as a linear combination of
$e(\nu)'$s with non-negative integer coefficients. That is
$$
{\text{\rm ch}}_3(\lambda)=\sum_{\nu \prec \lambda,~ \nu \in X(T)_+}b_{\lambda\nu}e(\nu)
$$
with $b_{\lambda\lambda}=1$. In this way, we get the second matrix $B$.
\par
For example, we assume that $G$ is the simple algebraic group
of type $A_5$ and $\lambda=(2,1,2,1,2)$.

It is easy to see that
$$
{\frak x}={\frak x}_{2\rho-\lambda}={\frak x}_{(01010)}=f_2f_1f_3f_2
f_4^{(2)}f_3^{(2)}f_2f_1f_5^{(2)}f_4^{(2)}f_3f_2.
$$
For $\nu=(3,0,1,2,2)$, we have $\lambda-\nu=(-1, 1, 1, -1,0)=\alpha_2+\alpha_3$.
First we compute each of the set $\SS_{\nu}=\{f_2f_3{\frak x}, f_{23}{\frak x}\}$.
Then we compute the rank of the set $\SS_{\nu}$, which is equal to $2$. So we have
$\dim L(2,1,2,1,2)_{(3,0,1,2,2)}=2$. For $\mu\!=\!(2,0,1,1,3)$, we have \!$\lambda-\mu\!=\!(0, 1, 1, 0,-1)\!=\!
\alpha_1\!+\!2\alpha_2\!+\!2\alpha_3\!+\!\alpha_4$. We compute each of the  set \!$\SS_{\mu}\!=\!\{
f_1f_2f_2f_3f_3f_4{\frak x}, f_1f_2f_2f_3 f_{34}{\frak x}, $
$f_1f_2f_3 f_{234}{\frak x},$
$ f_1f_2 f_{23}f_{34}{\frak x}, $
 $f_1 f_{23}f_{234}{\frak x},$
$~f_1f_2 f_{23}f_3f_4{\frak x},$
$ ~f_1 f_{23}f_{23}f_4{\frak x}, $
$~f_{12}f_2f_3f_3f_4{\frak x},$$ ~f_{12}f_2f_3 f_{34}{\frak x},$
$~f_{12}f_{23}f_{34}{\frak x},$
$ f_{12}f_{234}f_3{\frak x}, $ $ f_{12}f_{23}f_3f_4{\frak x},$
$ \!f_{123}f_{23}f_4{\frak x},
\!f_{123}f_2f_3f_4{\frak x}, $
$\!f_{123}f_2 f_{34}{\frak x}, \!f_{123}f_{234}{\frak x},$
$ \!f_{1234}f_{23}{\frak x},$
$ \!f_{1234}f_2f_3{\frak x} \!\}$,
 and then we compute the rank of the set $\SS_{\mu}$, which is equal to $13$. So we have
$\dim L(2,1,2,1,2)_{(2,0,1,1,3)}=13$. By this methods, we can calculate all multiplicity
$b_{\lambda\nu}.$   Finally,  we obtain the formal character of irreducible module
${\text{\rm ch}}_3(2,1,2,1,2).$
\par
When $\lambda$ lies in $X(T)_+$ but not in $X_1(T)$, we can also compute
the formal character ${\text{\rm ch}}_3(\lambda)$ by using the Steinberg tensor product
theorem. For $\lambda \in X(T)_+$, we have the unique decomposition
$$
\lambda=\lambda^0 + 3 \lambda^1 \qquad {\text {\rm with}}\quad \lambda^0 \in X_1(T),
\lambda^1 \in X(T)_+.
$$
Then the Steinberg tensor product theorem tells us that
$$
{\text{\rm ch}}_3(\lambda)={\text{\rm ch}}_3(\lambda^0) \cdot {\text{\rm ch}}_3(3 \lambda^1).
$$
Therefore, we can determine all characters ${\text{\rm ch}}_3 (\lambda)$ with $\lambda
\in X(T)_+$, provided that all characters ${\text{\rm ch}}_3 (\lambda)$ with $\lambda
\in X_1(T)$ are known. For example, when $\lambda = (0, 2, 0,0, 3)$, we have
$$\aligned
{\text{\rm ch}}_3(0, 2, 0, 0, 3)&={\text{\rm ch}}_3(0, 2, 0, 0,0)\cdot {\text{\rm ch}}_3(0, 0, 0, 0,3)\\
            &=(e(0, 2, 0, 0,0)+e(1, 0, 1, 0, 0)+ e(0, 0, 0, 1,0))\cdot e(0, 0, 0, 0,3)\\
            &=e(0, 2, 0, 0,3)+e(2, 0, 0,0,1)+e(1, 0,1 ,0, 3)+e(1,1, 0, 0, 2)\\
             \ \ &+e(0, 1, 0, 0,1)+e(0, 0, 0, 1,3)+e(0, 0,1, 0, 2).
\endaligned$$

Therefore, from the two matrices $A,B,$ we can easily get the third change
of basis matrix $D = AB^{-1}$ from $\{ {\text{\rm ch}}_3(\nu) \vert \nu \in X(T)_+
\}$ to $\{ {\text{\rm ch}}(\nu) \vert \nu \in X(T)_+ \}$, which is still a
triangular matrix with $1$ on its diagonal. The matrix $D$ gives the
decomposition patterns of various $H^0(\lambda)$ with $\lambda \in X(T)_+$.

We list the  matrix $D$ in the attached tables.
In all these tables, the left column indicates $\lambda$'s. For two weight  $\nu\prec\lambda \in
X(T)_+$, the number $d_{\lambda\nu}$ in tables is just the multiplicity of
composition factors $\left[ H^0(\lambda) : L(\nu)\right]$.

\section{FASTER ALGORITHM }

In our earlier paper \cite{YZ2001, YZ2002}, we compute the multiplicity $b_{\lambda\nu}$ one by one for a fixed
weight $\lambda.$
However,  noticing that some information  computing $b_{\lambda\nu}$ may
be useful to compute $b_{\lambda\mu}$ for $\nu\preccurlyeq \mu.$ So we
compute all possible $f_\beta,$ such that $\SS_{\lambda}=\{f_\beta{\frak x}\}$
spanning to the whole $L(\lambda),$ firstly.
Then we compute $\SS_{\lambda}=\{f_\beta{\frak x}\}$ in some ordering: if
$f_\beta=f_{\beta_1}f_{\beta_2},$ then we first obtain
$y_1=f_{\beta_2}{\frak x},$  save this result and compute
$y_2=f_\beta{\frak x}=f_{\beta_1}y_1$ instead of
computing $f_\beta{\frak x}=f_{\beta_1}f_{\beta_2}{\frak x},$
directly. In fact we only need compute $f_\beta y,$ for
some positive root $\beta $ and $y\in\SS_{\lambda}$  in one step.

For example, suppose  to compute $\{f_3f_4{\frak x},f_{23}f_4{\frak x}\},$ we
can compute  $y_1=f_4{\frak x}$ at the first step, and then compute
$y_2=f_3y_1,y_3=f_{23}y_1.$  In this way, we can avoid
much repeated work.

\section{WHY THESE RESULTS ARE TRUE}
In order to obtain the results the computer must work several days. So we
must be careful to avoid error.
There are facts to verity the results.

At firstly, we compute the dimension of weight space, then by Sternberg tensor
formula and Weyl formula we obtain the decomposition pattern of $H^0(\lambda).$
At last checking all the data  we find that

1). \,Symmetry of dimension of weight space.
Checking the results the two equations are satisfied:
\begin{eqnarray*}
&\dim L(\lambda_1,\lambda_2,\lambda_3,\lambda_2,\lambda_1)_{(\mu_1,\mu_2,\mu_3,\mu_4,\mu_5)}=
\dim L(\lambda_1,\lambda_2,\lambda_3,\lambda_2,\lambda_1)_{(\mu_5,\mu_4,\mu_3,\mu_2,\mu_1)},\\
&\dim L(\lambda_1,\lambda_2,\lambda_3,\lambda_4,\lambda_5)_{(\mu_1,\mu_2,\mu_3,\mu_4,\mu_5)}=
\dim L(\lambda_5,\lambda_4,\lambda_3,\lambda_2,\lambda_1)_{(\mu_5,\mu_4,\mu_3,\mu_2,\mu_1)}.
\end{eqnarray*}

2). \, Symmetry of composition factors.
From the $H^0(\lambda)'s$ decomposition  patterns,
the following equations are hold:
\begin{eqnarray*}
&&\left[ H^0(\lambda_1,\lambda_2,\lambda_3,\lambda_2,\lambda_1) : L(\mu_1,\mu_2,\mu_3,\mu_4,\mu_5)\right]\\
&&=\left[ H^0(\lambda_1,\lambda_2,\lambda_3,\lambda_2,\lambda_1) : L(\mu_5,\mu_4,\mu_3,\mu_2,\mu_1)\right].
\end{eqnarray*}

3). \, Positivity of multiplicity of composition factors. All the
multiplicity of composition factors we obtained  are nonnegative.


4). \, Linkage principle is hold. If the multiplicity of composition factors
 $\left[ H^0(\lambda) : L(\nu)\right]\neq 0, $ then
 we have $\mu\uparrow\lambda.$


From the representation theory of algebraic groups, all the above
results should be hold,  so the computational data is
compatible with the theory. Hence we may
accept these results.


\section{MAIN RESULTS}

{\bf Theorem }\quad
 When $G=SL(6,K)$, let $\Lambda=\{(2,2,2,2,2),$$(1,2,2,2,2),$ $(1,2,2,2,2),$ $
 (2,1,0,2,2), (2,2,0,1,2),(2,2,2,0,1),$
$(1,0,2,2,2),(2,0,1,2,2),$ $(2,2,1,0,2),$ $(0,2,2,2,2),$
$(2,2,2,2,0),(0,1,2,2,2),(2,2,2,1,0)$ $ \subset  X_1(T)$. Then
$H^0(\lambda)$ is an irreducible $G$-module for all $\lambda\in
\Lambda$ and the decomposition patterns of $H^0(\lambda)$ for all
$\lambda\in X_1(T) \backslash \Lambda$ are listed in Table 1-9.

{\bf Remark:\;} The table should be read as following. We list the weights in the first  collum
and write the multiplicity of composition factors as the others elements of tables.
For example, from the third row in table 1, we obtain $00200\ 0\ 1\ 1,$ this mean
$$
{\text{\rm ch}}(0,0,2,0,0)=0\cdot{\text{\rm ch}}_3(0,0,0,0,0)+1\cdot{\text{\rm ch}}_3(1,0,0,0,1)+1\cdot{\text{\rm ch}}_3(0,0,2,0,0).
$$
According to the symmetry of  $A_5,$ we need not list all results. For example,
we can obtain the decomposition pattern of $H^0(0,0,2,0,1)$ from
table 2:
$$
{\text{\rm ch}}(0,0,2,0,1)={\text{\rm ch}}_3(0,0,2,0,1)+{\text{\rm ch}}_3(1,0,0,0,2).
$$
So we also have
$$
{\text{\rm ch}}(1,0,2,0,0)={\text{\rm ch}}_3(1,0,2,0,0)+{\text{\rm ch}}_3(2,0,0,0,1).
$$


\center{\bf Table 1}
$$ \begin{smallmatrix}
00000& 1&   \cr
10001& 1&  1& \cr
00200& 0&  1&  1& \cr
01002& 0&  1&  0&  1& \cr
20010& 0&  1&  0&  0&  1& \cr
00111& 1&  1&  1&  1&  0&  1& \cr
11100& 1&  1&  1&  0&  1&  0&  1& \cr
00030& 1&  0&  0&  0&  0&  1&  0&  1& \cr
00103& 0&  0&  0&  1&  0&  1&  0&  0&  1& \cr
03000& 1&  0&  0&  0&  0&  0&  1&  0&  0&  1& \cr
30100& 0&  0&  0&  0&  1&  0&  1&  0&  0&  0&  1& \cr
11011& 2&  2&  1&  1&  1&  1&  1&  0&  0&  0&  0&  1& \cr
11003& 2&  1&  0&  1&  0&  1&  0&  0&  1&  0&  0&  1&  1& \cr
30011& 2&  1&  0&  0&  1&  0&  1&  0&  0&  0&  1&  1&  0&  1& \cr
00014& 0&  0&  1&  0&  0&  1&  0&  1&  1&  0&  0&  0&  0&  0&  1& \cr
41000& 0&  0&  1&  0&  0&  0&  1&  0&  0&  1&  1&  0&  0&  0&  0&  1& \cr
30003& 3&  1&  0&  0&  0&  0&  0&  0&  0&  0&  0&  1&  1&  1&  0&  0&  1& \cr
02020& 2&  0&  1&  0&  0&  1&  1&  1&  0&  1&  0&  1&  0&  0&  0&  0&  0&  1& \cr
10112& 3&  1&  1&  1&  0&  2&  1&  0&  1&  0&  0&  1&  1&  0&  0&  0&  0&  0&  1& \cr
21101& 3&  1&  1&  0&  1&  1&  2&  0&  0&  0&  1&  1&  0&  1&  0&  0&  0&  0&  0&  1& \cr
10031& 1&  0&  0&  0&  0&  1&  0&  1&  0&  0&  0&  0&  0&  0&  0&  0&  0&  0&  1&  0&  1& \cr
13001& 1&  0&  0&  0&  0&  0&  1&  0&  0&  1&  0&  0&  0&  0&  0&  0&  0&  0&  0&  1&  0&  1& \cr
00400& 0&  0&  1&  0&  0&  0&  0&  0&  0&  0&  0&  0&  0&  0&  0&  0&  0&  1&  0&  0&  0&  0&  1& \cr
02004& 2&  1&  1&  0&  0&  1&  0&  1&  1&  0&  0&  1&  1&  0&  1&  0&  0&  1&  0&  0&  0&  0&  0&  1& \cr
40020& 2&  1&  1&  0&  0&  0&  1&  0&  0&  1&  1&  1&  0&  1&  0&  1&  0&  1&  0&  0&  0&  0&  0&  0&  1& \cr
10023& 2&  0&  1&  0&  0&  1&  1&  1&  1&  0&  0&  0&  1&  0&  1&  0&  0&  0&  1&  0&  1&  0&  0&  0&  0&  1& \cr
32001& 2&  0&  1&  0&  0&  1&  1&  0&  0&  1&  1&  0&  0&  1&  0&  1&  0&  0&  0&  1&  0&  1&  0&  0&  0&  0&  1& \cr
01121& 2&  0&  1&  0&  0&  1&  1&  1&  0&  0&  0&  1&  1&  0&  0&  0&  0&  1&  1&  0&  1&  0&  0&  0&  0&  0&  0&  1& \cr
12110& 2&  0&  1&  0&  0&  1&  1&  0&  0&  1&  0&  1&  0&  1&  0&  0&  0&  1&  0&  1&  0&  1&  0&  0&  0&  0&  0&  0&  1& \cr
20202& 6&  1&  1&  0&  0&  1&  1&  0&  0&  0&  0&  1&  1&  1&  0&  0&  1&  0&  1&  1&  0&  0&  0&  0&  0&  0&  0&  0&  0&  1& \cr
00311& 0&  0&  1&  0&  0&  0&  1&  0&  0&  1&  0&  0&  0&  0&  0&  0&  0&  1&  0&  0&  0&  0&  1&  0&  0&  0&  0&  1&  0&  0&  1& \cr
11300& 0&  0&  1&  0&  0&  1&  0&  1&  0&  0&  0&  0&  0&  0&  0&  0&  0&  1&  0&  0&  0&  0&  1&  0&  0&  0&  0&  0&  1&  0&  0&  1& \cr
01113& 3&  1&  1&  0&  1&  1&  2&  1&  1&  1&  0&  1&  2&  0&  1&  0&  0&  1&  1&  0&  1&  0&  0&  1&  0&  1&  0&  1&  0&  0&  0&  0&  1& \cr
31110& 3&  1&  1&  1&  0&  2&  1&  1&  0&  1&  1&  1&  0&  2&  0&  1&  0&  1&  0&  1&  0&  1&  0&  0&  1&  0&  1&  0&  1&  0&  0&  0&  0&  1& \cr
40004& 3&  2&  0&  0&  0&  0&  0&  0&  0&  0&  0&  1&  1&  1&  0&  0&  1&  1&  0&  0&  0&  0&  0&  1&  1&  0&  0&  0&  0&  0&  0&  0&  0&  0&  1& \cr
00303& 0&  0&  0&  0&  1&  0&  1&  0&  0&  1&  0&  0&  0&  0&  0&  0&  0&  0&  0&  0&  0&  0&  0&  0&  0&  0&  0&  1&  0&  0&  1&  0&  1&  0&  0&  1& \cr
30300& 0&  0&  0&  1&  0&  1&  0&  1&  0&  0&  0&  0&  0&  0&  0&  0&  0&  0&  0&  0&  0&  0&  0&  0&  0&  0&  0&  0&  1&  0&  0&  1&  0&  1&  0&  0&  1& \cr
01032& 1&  0&  0&  0&  0&  0&  1&  0&  0&  0&  0&  0&  1&  0&  0&  0&  0&  0&  1&  0&  2&  0&  0&  0&  0&  1&  0&  1&  0&  0&  0&  0&  1&  0&  0&  0&  0&  1& \cr
23010& 1&  0&  0&  0&  0&  1&  0&  0&  0&  0&  0&  0&  0&  1&  0&  0&  0&  0&  0&  1&  0&  2&  0&  0&  0&  0&  1&  0&  1&  0&  0&  0&  0&  1&  0&  0&  0&  0&  1& \cr
11211& 6&  0&  2&  0&  0&  2&  2&  1&  0&  1&  0&  1&  1&  1&  0&  0&  0&  2&  1&  1&  1&  1&  1&  0&  0&  0&  0&  1&  1&  1&  1&  1&  0&  0&  0&  0&  0&  0&  0&  1& \cr
11203& 8&  1&  2&  0&  1&  1&  3&  0&  0&  2&  1&  1&  2&  1&  0&  0&  1&  1&  1&  1&  1&  1&  0&  1&  0&  1&  0&  1&  0&  1&  1&  0&  1&  0&  0&  1&  0&  0&  0&  1&  1& \cr
30211& 8&  1&  2&  1&  0&  3&  1&  2&  1&  0&  0&  1&  1&  2&  0&  0&  1&  1&  1&  1&  1&  1&  0&  0&  1&  0&  1&  0&  1&  1&  0&  1&  0&  1&  0&  0&  1&  0&  0&  1&  0&  1& \cr
30203&10&  2&  4&  0&  0&  2&  2&  1&  1&  1&  1&  1&  2&  2&  1&  1&  2&  1&  1&  1&  1&  1&  0&  1&  1&  1&  1&  0&  0&  1&  0&  0&  0&  0&  1&  0&  0&  0&  0&  1&  1&  1&  1& \cr
21212&10&  3&  7&  1&  1&  3&  3&  3&  1&  3&  1&  1&  2&  2&  1&  1&  1&  2&  1&  1&  2&  2&  1&  1&  1&  1&  1&  1&  1&  1&  1&  1&  1&  1&  0&  1&  1&  1&  1&  2&  1&  1&  1&  1& \cr
\end{smallmatrix} $$
\center{\bf Table 2}
$$
\begin{smallmatrix}
&\cr
00001& 1& \cr
12000& 1&  1& \cr
31000& 0&  1&  1& \cr
00120& 1&  0&  0&  1& \cr
11020& 2&  1&  0&  1&  1& \cr
00104& 0&  0&  0&  1&  0&  1& \cr
30020& 2&  1&  1&  0&  1&  0&  1& \cr
22001& 1&  1&  1&  0&  0&  0&  0&  1& \cr
11004& 2&  0&  0&  1&  1&  1&  0&  0&  1& \cr
10121& 1&  0&  0&  1&  1&  0&  0&  0&  0&  1& \cr
21110& 2&  1&  1&  1&  1&  0&  1&  1&  0&  0&  1& \cr
30004& 3&  0&  0&  0&  1&  0&  1&  0&  1&  0&  0&  1& \cr
20300& 0&  0&  0&  1&  0&  0&  0&  0&  0&  0&  1&  0&  1& \cr
10113& 2&  1&  0&  1&  1&  1&  0&  0&  1&  1&  0&  0&  0&  1& \cr
13010& 0&  0&  0&  0&  0&  0&  0&  1&  0&  0&  1&  0&  0&  0&  1& \cr
10032& 0&  0&  0&  0&  0&  0&  0&  0&  0&  1&  0&  0&  0&  1&  0&  1& \cr
20211& 2&  0&  0&  1&  1&  0&  1&  1&  0&  1&  1&  0&  1&  0&  0&  0&  1& \cr
01122& 1&  1&  0&  0&  1&  0&  0&  0&  1&  1&  0&  0&  0&  1&  0&  1&  0&  1& \cr
20203& 4&  1&  1&  0&  1&  0&  1&  1&  1&  1&  0&  1&  0&  1&  0&  0&  1&  0&  1& \cr
00312& 0&  1&  0&  0&  0&  0&  0&  0&  0&  0&  0&  0&  0&  0&  0&  0&  0&  1&  0&  1& \cr
11212& 4&  2&  1&  1&  1&  0&  1&  1&  1&  1&  1&  0&  1&  1&  1&  1&  1&  1&  1&  1&  1& \cr
\end{smallmatrix}
\
\begin{smallmatrix}
10002& 1& \cr
20100& 0&  1& \cr
00201& 1&  0&  1& \cr
20011& 1&  1&  0&  1& \cr
20003& 1&  0&  0&  1&  1& \cr
11101& 1&  1&  1&  1&  0&  1& \cr
03001& 0&  0&  0&  0&  0&  1&  1& \cr
30101& 0&  1&  0&  1&  0&  1&  0&  1& \cr
02110& 0&  0&  1&  1&  0&  1&  1&  0&  1& \cr
10202& 1&  0&  1&  1&  1&  1&  0&  0&  0&  1& \cr
41001& 0&  0&  1&  0&  0&  1&  1&  1&  0&  0&  1& \cr
01300& 0&  0&  1&  0&  0&  0&  0&  0&  1&  0&  0&  1& \cr
10040& 0&  0&  0&  0&  0&  0&  0&  0&  0&  1&  0&  0&  1& \cr
40110& 1&  0&  1&  1&  0&  1&  1&  1&  1&  0&  1&  0&  0&  1& \cr
01211& 0&  0&  1&  1&  0&  1&  1&  0&  1&  1&  0&  1&  0&  0&  1& \cr
01130& 0&  0&  0&  1&  1&  0&  0&  0&  0&  1&  0&  0&  1&  0&  1&  1& \cr
01203& 0&  1&  0&  1&  1&  1&  1&  0&  0&  1&  0&  0&  0&  0&  1&  0&  1& \cr
00320& 0&  0&  0&  0&  0&  0&  1&  0&  0&  0&  0&  1&  0&  0&  1&  1&  0&  1& \cr
01041& 0&  0&  0&  0&  1&  0&  0&  0&  0&  1&  0&  1&  1&  0&  1&  1&  1&  0&  1& \cr
11220& 0&  0&  1&  1&  1&  1&  1&  0&  1&  1&  0&  1&  1&  0&  1&  1&  0&  1&  0&  1& \cr
30220& 1&  0&  2&  1&  1&  1&  0&  0&  1&  1&  0&  0&  1&  1&  0&  0&  0&  0&  0&  1&  1& \cr
21221& 2&  1&  3&  1&  2&  1&  1&  0&  1&  1&  0&  1&  1&  1&  1&  1&  1&  1&  1&  1&  1&  1& \cr
\end{smallmatrix} $$
\center{\bf Table 3}
$$ \begin{smallmatrix}
10010& 1& \cr
01100& 1&  1& \cr
01011& 1&  1&  1& \cr
01003& 0&  0&  1&  1& \cr
50000& 0&  1&  0&  0&  1& \cr
00112& 0&  1&  1&  1&  0&  1& \cr
00031& 0&  0&  0&  0&  0&  1&  1& \cr
11012& 1&  1&  1&  1&  0&  1&  0&  1& \cr
00023& 0&  1&  0&  1&  0&  1&  1&  0&  1& \cr
30012& 1&  0&  0&  0&  0&  0&  0&  1&  0&  1& \cr
02021& 0&  1&  0&  0&  0&  1&  1&  1&  0&  0&  1& \cr
21102& 1&  1&  0&  0&  0&  1&  0&  1&  0&  1&  0&  1& \cr
02013& 1&  1&  0&  1&  0&  1&  1&  1&  1&  0&  1&  0&  1& \cr
12200& 0&  1&  0&  0&  0&  0&  0&  0&  0&  0&  0&  0&  0&  1& \cr
13002& 0&  0&  0&  0&  0&  0&  0&  0&  0&  0&  0&  1&  0&  0&  1& \cr
40021& 1&  1&  0&  0&  1&  0&  0&  1&  0&  1&  1&  0&  0&  0&  0&  1& \cr
31200& 0&  1&  1&  0&  1&  0&  0&  0&  0&  0&  0&  0&  0&  1&  0&  0&  1& \cr
32002& 0&  1&  0&  0&  1&  1&  0&  0&  0&  1&  0&  1&  0&  0&  1&  0&  0&  1& \cr
12111& 0&  2&  0&  0&  0&  1&  0&  1&  0&  1&  1&  1&  0&  1&  1&  0&  0&  0&  1& \cr
40013& 2&  0&  0&  0&  0&  0&  0&  1&  0&  1&  1&  0&  1&  0&  0&  1&  0&  0&  0&  1& \cr
23100& 0&  0&  0&  0&  0&  0&  0&  0&  0&  0&  0&  0&  0&  1&  0&  0&  1&  0&  0&  0&  1& \cr
12103& 1&  2&  0&  0&  0&  1&  0&  1&  0&  1&  1&  1&  1&  0&  1&  0&  0&  0&  1&  0&  0&  1& \cr
31111& 1&  3&  1&  1&  1&  2&  1&  1&  0&  2&  1&  1&  0&  1&  1&  1&  1&  1&  1&  0&  0&  0&  1& \cr
31103& 2&  4&  0&  1&  1&  2&  1&  1&  1&  2&  1&  1&  1&  0&  1&  1&  0&  1&  1&  1&  0&  1&  1&  1& \cr
23011& 0&  2&  0&  0&  1&  1&  0&  0&  0&  1&  0&  1&  0&  1&  2&  0&  1&  1&  1&  0&  1&  0&  1&  0&  1& \cr
23003& 0&  3&  0&  0&  1&  1&  0&  0&  0&  1&  0&  1&  0&  0&  2&  0&  0&  1&  1&  0&  0&  1&  1&  1&  1&  1& \cr
22112& 3&  8&  1&  1&  3&  2&  1&  1&  1&  2&  2&  1&  1&  1&  2&  1&  1&  1&  2&  1&  1&  1&  2&  1&  1&  1&  1& \cr
\end{smallmatrix} $$
\center{\bf Table 4}
$$ \begin{smallmatrix}
10210& 1& \cr
02221& 1&  1& \cr
\end{smallmatrix}
\
\begin{smallmatrix}
21021& 1& \cr
21013& 1&  1& \cr
12022& 1&  1&  1& \cr
\end{smallmatrix}
\
\begin{smallmatrix}
02102& 1& \cr
40102& 1&  1& \cr
22120& 1&  1&  1& \cr
\end{smallmatrix}
\
\begin{smallmatrix}
22010& 1& \cr
10122& 0&  1& \cr
20212& 1&  1&  1& \cr
\end{smallmatrix}
$$
\center{\bf Table 5}
$$ \begin{smallmatrix}
10012& 1& \cr
01102& 1&  1& \cr
50002& 0&  1&  1& \cr
02201& 0&  1&  0&  1& \cr
40201& 1&  1&  1&  1&  1& \cr
12210& 0&  1&  0&  1&  0&  1& \cr
24001& 0&  0&  0&  1&  1&  0&  1& \cr
31210& 1&  2&  0&  1&  1&  1&  0&  1& \cr
23110& 0&  2&  1&  1&  1&  1&  1&  1&  1& \cr
22300& 0&  1&  0&  0&  0&  1&  0&  1&  1&  1& \cr
22211& 1&  2&  1&  1&  1&  1&  1&  1&  1&  1&  1& \cr
\end{smallmatrix}
\
\begin{smallmatrix}
10100& 1& \cr
10011& 1&  1& \cr
10003& 0&  1&  1& \cr
01101& 1&  1&  0&  1& \cr
00202& 0&  1&  1&  1&  1& \cr
20012& 1&  1&  1&  0&  0&  1& \cr
50001& 0&  0&  0&  1&  0&  0&  1& \cr
00040& 0&  0&  0&  0&  1&  0&  0&  1& \cr
11102& 1&  1&  1&  1&  1&  1&  0&  0&  1& \cr
02200& 0&  0&  0&  1&  0&  0&  0&  0&  0&  1& \cr
03002& 0&  0&  0&  0&  0&  0&  0&  0&  1&  0&  1& \cr
30102& 1&  0&  0&  0&  0&  1&  0&  0&  1&  0&  0&  1& \cr
02111& 0&  0&  0&  1&  1&  1&  0&  0&  1&  1&  1&  0&  1& \cr
40200& 0&  1&  0&  1&  0&  0&  1&  0&  0&  1&  0&  0&  0&  1& \cr
41002& 0&  0&  0&  1&  1&  0&  1&  0&  1&  0&  1&  1&  0&  0&  1& \cr
02030& 0&  0&  1&  0&  1&  1&  0&  1&  0&  0&  0&  0&  1&  0&  0&  1& \cr
02103& 1&  0&  1&  0&  1&  1&  0&  0&  1&  0&  1&  0&  1&  0&  0&  0&  1& \cr
40111& 1&  1&  1&  1&  1&  1&  1&  0&  1&  1&  1&  1&  1&  1&  1&  0&  0&  1& \cr
24000& 0&  0&  0&  0&  0&  0&  0&  0&  0&  1&  0&  0&  0&  1&  0&  0&  0&  0&  1& \cr
40030& 1&  0&  1&  0&  0&  1&  0&  0&  0&  0&  0&  0&  1&  0&  0&  1&  0&  1&  0&  1& \cr
40103& 2&  0&  1&  0&  1&  1&  0&  0&  1&  0&  1&  1&  1&  0&  1&  0&  1&  1&  0&  0&  1& \cr
12120& 0&  0&  1&  1&  1&  1&  0&  0&  1&  1&  1&  1&  1&  0&  0&  1&  0&  0&  0&  0&  0&  1& \cr
31120& 1&  1&  2&  2&  2&  1&  0&  1&  1&  1&  1&  1&  1&  1&  1&  1&  0&  1&  0&  1&  0&  1&  1& \cr
23020& 0&  0&  0&  2&  1&  0&  1&  0&  1&  1&  2&  1&  0&  1&  1&  0&  0&  0&  1&  0&  0&  1&  1&  1& \cr
22121& 3&  2&  3&  3&  2&  1&  1&  1&  1&  1&  2&  1&  2&  1&  1&  1&  1&  1&  0&  1&  1&  1&  1&  1&  1& \cr
\end{smallmatrix} $$
\center{\bf Table 6}
$$
\begin{smallmatrix}
12010& 1& \cr
10212& 1&  1& \cr
10131& 0&  1&  1& \cr
20221& 1&  1&  1&  1& \cr
\end{smallmatrix}
\
\begin{smallmatrix}
02101& 1& \cr
01202& 1&  1& \cr
01040& 0&  1&  1& \cr
21220& 1&  1&  1&  1& \cr
\end{smallmatrix}
\
\begin{smallmatrix}
01012& 1& \cr
12201& 0&  1& \cr
31201& 1&  1&  1& \cr
23101& 0&  1&  1&  1& \cr
22202& 1&  1&  1&  1&  1& \cr
\end{smallmatrix}
\
\begin{smallmatrix}
20101& 1& \cr
01220& 0&  1& \cr
01204& 1&  1&  1& \cr
01042& 0&  1&  1&  1& \cr
21222& 1&  1&  1&  1&  1& \cr
\end{smallmatrix}
\
\begin{smallmatrix}
20002& 1& \cr
10201& 1&  1& \cr
01210& 0&  1&  1& \cr
02220& 0&  1&  1&  1& \cr
12221& 1&  1&  0&  1&  1& \cr
\end{smallmatrix}
$$
\center{\bf Table 7}
$$ \begin{smallmatrix}
&\cr
&\cr
00002& 1& \cr
21000& 0&  1& \cr
00210& 1&  0&  1& \cr
20020& 1&  1&  0&  1& \cr
12001& 1&  1&  0&  0&  1& \cr
31001& 0&  1&  0&  0&  1&  1& \cr
11110& 1&  1&  1&  1&  1&  0&  1& \cr
20004& 1&  0&  0&  1&  0&  0&  0&  1& \cr
10300& 0&  0&  1&  0&  0&  0&  1&  0&  1& \cr
03010& 0&  0&  0&  0&  1&  0&  1&  0&  0&  1& \cr
30110& 1&  1&  0&  1&  1&  1&  1&  0&  0&  0&  1& \cr
10211& 1&  0&  1&  1&  1&  0&  1&  0&  1&  0&  0&  1& \cr
10130& 0&  0&  0&  1&  0&  0&  0&  0&  0&  0&  0&  1&  1& \cr
10203& 1&  1&  0&  1&  1&  0&  0&  1&  0&  0&  0&  1&  0&  1& \cr
10041& 0&  0&  0&  0&  0&  0&  0&  0&  1&  0&  0&  1&  1&  1&  1& \cr
01212& 0&  1&  0&  1&  1&  0&  1&  0&  1&  1&  0&  1&  0&  1&  0&  1& \cr
20220& 1&  0&  0&  1&  1&  0&  1&  0&  1&  0&  1&  1&  1&  0&  0&  0&  1& \cr
01131& 0&  0&  0&  1&  0&  0&  0&  1&  1&  0&  0&  1&  1&  1&  1&  1&  0&  1& \cr
00321& 0&  0&  0&  0&  0&  0&  0&  0&  1&  1&  0&  0&  0&  0&  0&  1&  0&  1&  1& \cr
11221& 2&  1&  1&  1&  1&  0&  1&  1&  2&  1&  1&  1&  1&  1&  1&  1&  1&  1&  1&  1&
\end{smallmatrix}
\;\;\;
\begin{smallmatrix}
00010& 1& \cr
02000& 1&  1& \cr
40000& 0&  1&  1& \cr
01020& 1&  1&  0&  1& \cr
01004& 0&  0&  0&  1&  1& \cr
00121& 0&  0&  0&  1&  0&  1& \cr
00113& 0&  1&  0&  1&  1&  1&  1& \cr
11021& 1&  1&  0&  1&  0&  1&  0&  1& \cr
00032& 0&  0&  0&  0&  0&  1&  1&  0&  1& \cr
11013& 2&  1&  0&  1&  1&  1&  1&  1&  0&  1& \cr
30021& 2&  1&  1&  0&  0&  0&  0&  1&  0&  0&  1& \cr
21200& 0&  1&  1&  1&  0&  0&  0&  0&  0&  0&  0&  1& \cr
22002& 1&  1&  1&  0&  0&  0&  0&  0&  0&  0&  0&  0&  1& \cr
30013& 3&  0&  0&  0&  0&  0&  0&  1&  0&  1&  1&  0&  0&  1& \cr
13100& 0&  0&  0&  0&  0&  0&  0&  0&  0&  0&  0&  1&  0&  0&  1& \cr
21111& 2&  1&  1&  1&  0&  1&  0&  1&  0&  0&  1&  1&  1&  0&  0&  1& \cr
02022& 1&  1&  0&  0&  0&  1&  1&  1&  1&  1&  0&  0&  0&  0&  0&  0&  1& \cr
21103& 4&  1&  1&  0&  0&  1&  1&  1&  0&  1&  1&  0&  1&  1&  0&  1&  0&  1& \cr
13011& 0&  0&  1&  0&  0&  0&  0&  0&  0&  0&  0&  1&  1&  0&  1&  1&  0&  0&  1& \cr
13003& 1&  0&  1&  0&  0&  0&  0&  0&  0&  0&  0&  0&  1&  0&  0&  1&  0&  1&  1&  1& \cr
12112& 4&  3&  3&  1&  0&  1&  1&  1&  0&  1&  1&  1&  1&  1&  1&  2&  1&  1&  1&  1&  1&
\end{smallmatrix} $$
\center{\bf Table 8}
$$ \begin{smallmatrix}
00100& 1& \cr
00011& 1&  1& \cr
11000& 1&  0&  1& \cr
00003& 0&  1&  0&  1& \cr
30000& 0&  0&  1&  0&  1& \cr
10020& 1&  1&  1&  0&  0&  1& \cr
02001& 1&  1&  1&  0&  0&  0&  1& \cr
01110& 1&  1&  1&  0&  0&  1&  1&  1& \cr
10004& 0&  1&  0&  1&  0&  1&  0&  0&  1& \cr
40001& 0&  0&  1&  0&  1&  0&  1&  0&  0&  1& \cr
00300& 0&  0&  0&  0&  0&  0&  0&  1&  0&  0&  1& \cr
00211& 0&  1&  0&  1&  0&  1&  1&  1&  0&  0&  1&  1& \cr
11200& 0&  0&  1&  0&  1&  1&  1&  1&  0&  0&  1&  0&  1& \cr
20021& 1&  1&  1&  1&  1&  1&  0&  0&  0&  0&  0&  0&  0&  1& \cr
12002& 1&  1&  1&  1&  1&  0&  1&  0&  0&  0&  0&  0&  0&  0&  1& \cr
00130& 0&  0&  0&  0&  0&  1&  0&  0&  0&  0&  0&  1&  0&  0&  0&  1& \cr
00203& 0&  1&  1&  1&  0&  1&  1&  0&  1&  0&  0&  1&  0&  0&  0&  0&  1& \cr
03100& 0&  0&  0&  0&  0&  0&  1&  0&  0&  0&  0&  0&  1&  0&  0&  0&  0&  1& \cr
30200& 0&  1&  1&  0&  1&  1&  1&  0&  0&  1&  0&  0&  1&  0&  0&  0&  0&  0&  1& \cr
20013& 1&  1&  0&  1&  0&  1&  0&  0&  1&  0&  0&  0&  0&  1&  0&  0&  0&  0&  0&  1& \cr
31002& 1&  0&  1&  0&  1&  0&  1&  0&  0&  1&  0&  0&  0&  0&  1&  0&  0&  0&  0&  0&  1& \cr
11111& 1&  1&  1&  1&  1&  2&  2&  1&  0&  0&  1&  1&  1&  1&  1&  0&  0&  0&  0&  0&  0&  1& \cr
11030& 0&  0&  0&  1&  0&  1&  0&  0&  0&  0&  0&  1&  0&  1&  0&  1&  0&  0&  0&  0&  0&  1&  1& \cr
03011& 0&  0&  0&  0&  1&  0&  1&  0&  0&  0&  0&  0&  1&  0&  1&  0&  0&  1&  0&  0&  0&  1&  0&  1& \cr
11103& 2&  1&  1&  2&  1&  1&  1&  0&  1&  0&  0&  1&  0&  1&  1&  0&  1&  0&  0&  1&  0&  1&  0&  0&  1& \cr
30111& 2&  1&  1&  1&  2&  1&  1&  0&  0&  1&  0&  0&  1&  1&  1&  0&  0&  0&  1&  0&  1&  1&  0&  0&  0&  1& \cr
00041& 0&  0&  0&  0&  0&  0&  0&  0&  0&  0&  1&  1&  0&  0&  0&  1&  1&  0&  0&  0&  0&  0&  0&  0&  0&  0&  1& \cr
14000& 0&  0&  0&  0&  0&  0&  0&  0&  0&  0&  1&  0&  1&  0&  0&  0&  0&  1&  1&  0&  0&  0&  0&  0&  0&  0&  0&  1& \cr
30030& 1&  0&  0&  1&  1&  0&  0&  0&  0&  0&  0&  0&  0&  1&  0&  0&  0&  0&  0&  0&  0&  1&  1&  0&  0&  1&  0&  0&  1& \cr
03003& 1&  0&  0&  1&  1&  0&  0&  0&  0&  0&  0&  0&  0&  0&  1&  0&  0&  0&  0&  0&  0&  1&  0&  1&  1&  0&  0&  0&  0&  1& \cr
30103& 3&  0&  0&  1&  1&  0&  0&  0&  0&  0&  0&  0&  0&  1&  1&  0&  0&  0&  0&  1&  1&  1&  0&  0&  1&  1&  0&  0&  0&  0&  1& \cr
02112& 1&  0&  1&  1&  3&  1&  1&  0&  0&  0&  1&  1&  1&  1&  1&  0&  1&  1&  0&  1&  0&  2&  0&  1&  1&  0&  0&  0&  0&  1&  0&  1& \cr
21120& 1&  1&  0&  3&  1&  1&  1&  0&  0&  0&  1&  1&  1&  1&  1&  1&  0&  0&  1&  0&  1&  2&  1&  0&  0&  1&  0&  0&  1&  0&  0&  0&  1& \cr
02031& 0&  0&  0&  1&  1&  1&  0&  0&  1&  0&  1&  1&  0&  1&  0&  1&  1&  0&  0&  1&  0&  1&  1&  0&  0&  0&  1&  0&  0&  0&  0&  1&  0&  1& \cr
13020& 0&  0&  0&  1&  1&  0&  1&  0&  0&  1&  1&  0&  1&  0&  1&  0&  0&  1&  1&  0&  1&  1&  0&  1&  0&  0&  0&  1&  0&  0&  0&  0&  1&  0&  1& \cr
12121& 3&  2&  2&  4&  4&  2&  2&  1&  1&  1&  2&  1&  1&  1&  1&  1&  1&  1&  1&  1&  1&  3&  1&  1&  1&  1&  0&  0&  1&  1&  1&  1&  1&  1&  1&  1& \cr
\end{smallmatrix} $$
\center{\bf Table 9}
$$
\begin{smallmatrix}
00122& 1& \cr
22100& 0&  1& \cr
11022& 1&  0&  1& \cr
22011& 0&  1&  0&  1& \cr
30022& 0&  0&  1&  0&  1& \cr
22003& 0&  0&  0&  1&  0&  1& \cr
21112& 1&  1&  1&  1&  1&  1&  1& \cr
\end{smallmatrix}
\
\begin{smallmatrix}
01010& 1& \cr
21012& 1&  1& \cr
12021& 0&  1&  1& \cr
12013& 1&  1&  1&  1& \cr
31021& 1&  1&  1&  0&  1& \cr
31013& 2&  1&  1&  1&  1&  1& \cr
22022& 3&  1&  2&  1&  1&  1&  1& \cr
\end{smallmatrix}
\
\begin{smallmatrix}
10101& 1& \cr
20102& 1&  1& \cr
02120& 0&  1&  1& \cr
02104& 1&  1&  1&  1& \cr
40120& 1&  1&  1&  0&  1& \cr
40104& 2&  1&  1&  1&  1&  1& \cr
22122& 3&  1&  2&  1&  1&  1&  1& \cr
\end{smallmatrix}
\
\begin{smallmatrix}
00022& 1& \cr
02012& 1&  1& \cr
40012& 0&  1&  1& \cr
12102& 0&  1&  0&  1& \cr
31102& 1&  1&  1&  1&  1& \cr
22200& 0&  0&  0&  0&  0&  1& \cr
23002& 0&  0&  0&  1&  1&  0&  1& \cr
22111& 1&  1&  1&  1&  1&  1&  1&  1& \cr
\end{smallmatrix}
$$

\medskip

\begin{flushleft} \noindent{\bf Acknowledgement:} This work was
supported  by the Natural Science Fund of Hohai
University(2084/409277,2084/407188) and the Fundamental Research
Funds for the Central Universities 2009B26914 and 2010B09714. The
authors wishes to thank Prof. Ye Jiachen  for his helpful advice.
\end{flushleft}








\end{document}